\documentclass[11pt,twoside,reqno]{amsart}
\usepackage{amsmath,amsfonts,amsthm,amssymb}
\usepackage{graphicx,psfrag,overpic}
\usepackage{bm}
\usepackage[all,cmtip,line]{xy}
\usepackage{array,setspace}
\usepackage{color}
\numberwithin{equation}{section}

\usepackage{graphics}
\usepackage{epsfig}
\newtheorem{claim}{\bf \t}[part]

\newtheorem{theorem}{Theorem}[section]

\newtheorem{remark}{Remark}[section]

\def\v{\varepsilon}

\def\t{\theta}

\def\r{\rho}

\newcommand \R{\mathbb{R}}

\begin{document}

\title[Subsonic-Sonic Limit to M-D Steady Euler Equations]
{Subsonic-Sonic Limit of Approximate Solutions to Multidimensional Steady Euler Equations}

\author{Gui-Qiang Chen}
\address{G.-Q. Chen,
School of Mathematical Sciences, Fudan University\\
Shanghai 200433, China;
Mathematical Institute, University of Oxford\\
Radcliffe Observatory Quarter, Woodstock Road, Oxford, OX2 6GG, UK;
Academy of Mathematics and Systems Science\\
Academia Sinica, Beijing, 100190, P. R. China}
\email{chengq@maths.ox.ac.uk}

\author{Fei-Min Huang}
\address{F.-M. Huang, Institute of Applied Mathematics\\
 Academy of Mathematics and Systems Science\\
Academia Sinica, Beijing, 100190, P. R. China}
\email{fhuang@amt.ac.cn}

\author{Tian-Yi Wang}
\address{T.-Y. Wang, Department of Mathematics, School of Science,
Wuhan University of Technology, Wuhan, Hubei, 430070, P. R. China;
Academy of Mathematics and Systems Science, Academia Sinica, Beijing, 100190, P. R. China;
Mathematical Institute, University of Oxford, Radcliffe Observatory Quarter,
Woodstock Road, Oxford, OX2 6GG, UK}
\email{tianyiwang@whut.edu.cn; wangtianyi@amss.ac.cn}
\date{\today}

\begin{abstract}
A compactness framework is established for approximate solutions to subsonic-sonic flows
governed by the steady full Euler equations for compressible fluids in arbitrary dimension.
The existing compactness frameworks for the two-dimensional irrotational case do not
directly apply for the steady full Euler
equations in higher dimensions.
The new compactness framework we develop
applies for both non-homentropic and rotational flows.
One of our main observations is that the compactness can be achieved
by using only natural weak estimates for the mass balance and
the vorticity,
along with the Bernoulli law and the entropy relation,
through a more delicate analysis on the phase space.
As direct applications, we establish two existence theorems for
multidimensional subsonic-sonic full Euler flows through
infinitely long nozzles.
\end{abstract}

\keywords{Multidimensional, subsonic-sonic limit, steady flow, full Euler equations, homentropic,
rotation, compactness framework, strong convergence, exact solutions, approximate solutions}
\subjclass[2010]{
35Q31; 
35M30; 
35L65; 
76N10; 
76G25; 
35B40; 
35D30
}

\maketitle

\section{Introduction}
The full Euler equations for steady compressible flows in $\R^n$ read
\begin{eqnarray}\label{1.5}
\begin{cases}
\mbox{div}\,(\rho u)=0,\\
\mbox{div}\,(\rho u\otimes u+ p I)=0,\\
\mbox{div}\, (\rho u E+ u p)=0,
\end{cases}
\end{eqnarray}
where $x=(x_1, \cdots, x_n)\in \R^n, n\ge 2$, $u=(u_1,\cdots,u_n)\in \R^n$ is the fluid velocity,
and
$$
q:=|u|=\Big(\sum_{i=1}^n u_i^2\Big)^{1/2}
$$
is the speed, while $\rho$, $p$, and $E$ represent the density,
pressure, and total energy,
respectively, and $I$ is the $n\times n$ unit matrix.
The nonnegative quantities $\rho$, $q$, $p$, and $E$ are not independent.

For ideal polytropic gas,
$$
E=\frac{q^2}{2}+\frac{p}{(\gamma-1)\rho}
$$
with adiabatic exponent $\gamma>1$.
In this case, the Bernoulli law is written as
\begin{equation}\label{1.6}
\frac{q^2}{2}+ h(\rho, p)=B,
\end{equation}
where $h(\rho, p)=\frac{\gamma p}{(\gamma-1)\rho}$ is the enthalpy,
and $B$ is a Bernoulli function determined by appropriate additional conditions (such as
boundary conditions and/or asymptotic conditions at infinity).
The sound speed of the flow is
\begin{equation}\label{1.7}
c=\sqrt{\frac{\gamma p}{\rho}},
\end{equation}
and the Mach number is  defined as
\begin{equation}\label{1.8}
M=\frac{q}{c}.
\end{equation}
Then, for a fixed Bernoulli function $B$,
there is a critical speed $q_{\rm cr}=\sqrt{2\frac{\gamma-1}{\gamma+1}B}$
such that, when $q\le q_{\rm cr}$, the flow is subsonic-sonic (that is, $M \le 1$);
otherwise, it is supersonic (that is, $M>1$).

\smallskip
When the flow is homentropic, the pressure is a function of the density, and $(\ref{1.5})$
then reduces to
\begin{eqnarray}\label{1.1}
\begin{cases}
\mbox{div}(\rho u)=0,\\[2mm]
\mbox{div}(\rho u\otimes u+pI)=0.
\end{cases}
\end{eqnarray}
As usual, we require
$$
p'(\rho)>0,  \quad 2p'(\rho) + \rho p''(\rho) > 0 \qquad \mbox{for $\rho>0$},
$$
which include the $\gamma$-law flows with  $p=\kappa \rho^\gamma$, $\gamma>1$ and $\kappa>0$,
and the isothermal flows with $p=\kappa \rho$; see \cite{Courant-Friedrichs}.
In this case, the Bernoulli law has the same formula as $(\ref{1.6})$,
while the enthalpy $h(p(\rho),\rho)=h(\rho)$ with $h'(\rho)=\frac{p'(\rho)}{\rho}$.
The sound speed of the flow is $c=\sqrt{p'(\rho)}$,
and the Mach number is  defined as $M=\frac{q}{c}$.
Then, for a fixed Bernoulli function $B$,
the condition $2p'(\rho) + \rho p''(\rho) > 0$ for $\rho>0$ implies the existence of
a unique critical speed $q_{\rm cr}=q_{\rm cr}(B)$ so that
the flow can be classified by subsonic-sonic or supersonic
by $q\leq q_{\rm cr}$ or $q>q_{\rm cr}$, respectively.

\smallskip
It is well known that the steady Euler equations for compressible fluids are of composite-mixed
type, which is determined by the Mach number $M$. That is, the system can be reduced
to a system such that two of the equations are elliptic-hyperbolic mixed: elliptic when $M < 1$
and hyperbolic when $M > 1$, while the other $n$ equations are hyperbolic.
For the homentropic case, there are two mixed characteristics and $n-1$ hyperbolic characteristics.

\smallskip
During the 1950s, the effort on system \eqref{1.1} was focused mainly on the irrotational case,
namely when $u$ is constrained to satisfy the additional equation $\mbox{curl} \,u=0$.
Since the equations of uniform subsonic flow possess ellipticity,
solutions have better regularity
than those corresponding to transonic or supersonic flow.
The airfoil problem for irrotational, two-dimensional subsonic flow was solved; {\it cf.}
Frankl-Keldysh \cite{Frankl},  Shiffman \cite{Shiffman2},  Bers \cite{Bers2}, and Finn-Gilbarg \cite{Finn1}.
The first result for three-dimensional subsonic flow past
an obstacle was given by Finn
and Gilbarg \cite{Gilbarg1}
under some restrictions on the Mach number.
Dong \cite{Dong1} and Dong-Ou \cite{Dong2} extended the results to the maximum Mach number $M < 1$ for
the arbitrarily dimensional case.
Also see Du-Xin-Yan \cite{Xin3} for the construction of a
smooth uniform subsonic flow in an infinitely long nozzle
in $\R^n, n\geq 2$.

For the rotational case,
the global existence of homentropic subsonic flow through
two-dimensional infinitely long nozzles was proved in Xie-Xin \cite{Xin4}.
The result was also extended to the two-dimensional periodic nozzles in Chen-Xie \cite{CX}
and to axisymmetric nozzles in Du-Duan \cite{DD}.
For the full Euler flow, the first result was given by Chen-Deng-Xiang \cite{ChenDX}
for two-dimensional infinitely long nozzles.
Bae \cite{Bea} showed the stability of contact discontinuities for subsonic full Euler flow
in the two-dimensional, infinitely long nozzles.
Duan-Luo \cite{Duan-Luo} recently considered the axisymmetric nozzle
problem for the smooth subsonic flow.

\smallskip
On the other hand, few results are currently known for the cases of
subsonic-sonic flow
and transonic flow, since the uniform ellipticity is lost and
shocks may be present. That is, smooth solutions may not exist.
Instead, one must consider weak solutions.
Morawetz \cite{Morawetz2,Morawetz3} introduced an approach via
compensated compactness to analyze irrotational steady flow
of the Euler
equations. Indeed, Morawetz established a compactness framework
under the assumption that the solutions are free of stagnation points
and cavitation points. Morawetz's result has been improved
by Chen-Slemrod-Wang \cite{Chen8}, who have shown that the approximate
solutions away from cavitation are constructed
by a viscous perturbation.

The compactness framework for subsonic-sonic
irrotational flow allowing for stagnation in two dimensions
was due to
Chen-Dafermos-Slemrod-Wang \cite{Chen6}
by combining the mass conservation, momentum, and
irrotational equations.
The key observation in \cite{Chen6}
is that the two-dimensional steady
flow can be regarded as a one-dimensional unsteady system of conservation
laws, that is, one of the spatial variables can be regarded as the time variable,
so that the div-curl lemma can be applied to the two momentum equations. In fact,
the momentum equations are first employed in \cite{Chen6} to reduce
the support of the corresponding Young measure to two points,
and then the irrotational equation and the mass equation are used
to reduce the Young measure to a Dirac measure.
Using a similar idea, Xie-Xin \cite{Xin1} investigated the subsonic-sonic limit
of the two-dimensional irrotational, infinitely long nozzle problem.
Later, in  \cite{Xin2}, they extended the result to the three-dimensional
axisymmetric flow
through an axisymmetric nozzle.
The compactness framework in the multidimensional irrotational case was
established in Huang-Wang-Wang \cite{Huang-Wang-Wang}.

\smallskip
The compactness framework established for irrotational flow
no longer applies directly for the steady full Euler equations
in $\R^n$ with $n\geq 2$. When $n\ge 3$, the equations cannot be reduced
to a one-dimensional system of
conservation laws. More importantly, the div-curl lemma is no longer valid
for the momentum equations, due to the presence of linear characteristics.
One of our main observations is
that it is still possible to achieve the same compactness result, {\it i.e.},
to reduce the Young measure to a Dirac
measure, by using only natural weak estimates for the mass balance and the vorticity,
along with the Bernoulli law and entropy relation, through
a more delicate analysis on the phase space.
In particular,
the Bernoulli function and entropy function play a key role in our proof.
We then establish a compactness framework for approximate solutions
for steady full Euler flows in
arbitrary dimension.

\smallskip
The rest of this paper is organized as follows.
In Section 2, we establish the compactness framework for subsonic-sonic approximate
solutions to subsonic-sonic flows governed by the steady full Euler equations
for compressible fluids in $\R^n$ with $n\geq2$,
as well as by the steady homentropic Euler equations with weaker conditions.
In Sections 3--4, we give two direct applications of the compactness framework
to establish the existence of subsonic-sonic full Euler flow
through infinitely long nozzles
in $\R^n$ with $n\geq 2$.

\section{Compactness Framework for Approximate Steady Full Euler Flows}

In this section, we establish the compensated compactness framework for approximate solutions
of the steady full Euler equations in $\R^n$ with $n\ge 2$ with the form:
\begin{eqnarray}\label{3.1}
\begin{cases}
\mbox{div}(\rho^\varepsilon u^\varepsilon)=e_1(\varepsilon),\\[1mm]
\mbox{div}(\rho^\varepsilon u^\varepsilon\otimes u^\varepsilon+ p^\varepsilon I)=e_2(\varepsilon),\\[1mm]
\mbox{div} (\rho^\varepsilon u^\varepsilon E^\varepsilon+ u^\varepsilon p^\varepsilon)=e_3(\varepsilon),
\end{cases}
\end{eqnarray}
where $e_1(\varepsilon)$, $e_2(\varepsilon)=(e_{21}(\varepsilon), \cdots, e_{2n}(\varepsilon))^\top$,
and $e_3(\varepsilon)$ are sequences of functions depending on the parameter $\varepsilon$.

Let a sequence of functions $\rho^\varepsilon(x)$, $u^\varepsilon(x)=(u^\varepsilon_1,
\cdots, u^\varepsilon_n)(x)$, and $p^\varepsilon(x)$ be defined on an open subset $\Omega\subset
\mathbb{R}^n$ such that the following qualities:
\begin{eqnarray}
&&q^\varepsilon := |u^\varepsilon|=\sqrt{\sum_{i=1}^n (u_i^\varepsilon)^2}, \quad
c^\varepsilon := \sqrt{\frac{\gamma p^\varepsilon}{\rho^\varepsilon}}, \quad
M^\varepsilon := \frac{q^\varepsilon}{c^\varepsilon},\label{2.2a}\\[3mm]
&&B^\varepsilon := \frac{(q^\varepsilon)^2}{2}+\frac{\gamma p^\varepsilon}{(\gamma-1)\rho^\varepsilon}, \qquad
S^\varepsilon := \frac{\gamma p^\varepsilon}{(\gamma-1)(\rho^\varepsilon)^\gamma} \label{2.2b}
\end{eqnarray}
can be well defined and satisfy the following conditions:

\medskip
(A.1). $M^\varepsilon\leq 1$  {\it a.e.} in $\Omega$;

\vspace{2mm}

(A.2). $S^\varepsilon$ and $B^\varepsilon$ are uniformly bounded and,
for any compact set $K$, there exists a uniform constant $c(K)$ such
that $\inf\limits_{x\in K} S^\varepsilon(x)\ge c(K)>0$.
Moreover,  $(S^\varepsilon, B^\varepsilon)\to (\overline{S}, \overline{B})$ {\it a.e.} in $\Omega$;

\vspace{2mm}
(A.3). $\mbox{curl}\  u^\varepsilon$ and $e_1(\varepsilon)$ are in a compact set in $W_{loc}^{-1, p}$ for some $1<p \le 2$.

\medskip
\noindent
Then we have

\begin{theorem}[Compensated  compactness framework for the full Euler case] \label{thm3.1}
Let a sequence of functions $\rho^\varepsilon(x)$, $u^\varepsilon(x)=(u^\varepsilon_1,
\cdots, u^\varepsilon_n)(x)$, and $p^\varepsilon(x)$ satisfy conditions {\rm (A.1)}--{\rm (A.3)}.
Then there exists a subsequence (still labeled) $(\rho^\varepsilon, u^\varepsilon, p^\varepsilon)(x)$
such that
$$
\rho^\varepsilon(x)\to \rho(x), \quad u^\varepsilon(x)\to (u_1,
\cdots, u_n)(x), \quad p^\varepsilon(x)\to p(x)
\qquad \mbox{a.e. in $x \in \Omega$ as $\varepsilon\rightarrow 0$},
$$
and
$$
M(x):=\frac{q(x)}{c(x)}\leq 1  \qquad \mbox{a.e. $x\in \Omega$}.
$$
\end{theorem}

\noindent\textbf{Proof}.
We divide the proof into three steps.

\medskip
{\it Step 1. The strong convergence of $(\rho^\v, p^\v)$ follows from the
strong convergence of $q^\v=|u^\v|$}.

We employ \eqref{2.2b} to obtain
\begin{equation}\label{3.3f}
\frac{(q^\v)^2}{2}+S^\v(\rho^\v)^{\gamma-1}=B^\v.
\end{equation}

From these, the three variables $\rho^\v$, $p^\v$, and $q^\v$ are determined by one of them.
In other words, the pressure $p^\v$ and density $\rho^\v$ can be regarded as functions
of $q^\v$ through $(B^\v, S^\v)$ with
\begin{eqnarray}
&&\rho^\varepsilon=\rho(q^\varepsilon; B^\varepsilon, S^\varepsilon)
=\left(\frac{2B^\varepsilon- (q^\varepsilon)^2}{2S^\varepsilon}\right)^{\frac{1}{\gamma-1}},\label{3.3a1}\\[2mm]
&&p^\varepsilon=p(q^\varepsilon; B^\varepsilon, S^\varepsilon)
=\frac{\gamma-1}{2\gamma}
\frac{\left(2B^\varepsilon-(q^\varepsilon)^2\right)^{\frac{\gamma}{\gamma-1}}}{(2S^\varepsilon)^{\frac{1}{\gamma-1}}}.
\label{3.3a2}
\end{eqnarray}
Since $(B^\v, S^\v)$ strongly converge to $(\overline{B}, \overline{S})$  {\it a.e.},
the strong convergence of the density $\rho^\v$ and pressure $p^\v$
becomes a nature consequence of the strong convergence of the speed $q^\v$.

\medskip
{\it Step 2. The $H^{-1}_{loc}$--compactness}.
From the uniform boundedness of the Bernoulli function $B^\varepsilon$ and
the subsonic-sonic condition $M^\varepsilon\le1$,
it is easy to see
$$
q^\varepsilon\le \sup\limits_{\varepsilon}\sqrt{2\frac{\gamma-1}{\gamma+1}B^\varepsilon},
$$
which implies the uniform boundedness of
the velocity $u^\varepsilon(x)=(u^\varepsilon_1, \cdots, u^\varepsilon_n)(x)$.
This shows that
$$
(\mbox{curl}\,u^\varepsilon)_{ij}=\partial_{x_i} u^\varepsilon_j - \partial_{x_j} u^\varepsilon_i
\qquad\mbox{is bounded in $W^{-1, \infty}$}.
$$
On the other hand, $\mbox{curl}\, u^\varepsilon$ is compact in $W^{-1, p}$ for some $1< p \le 2$.
By the interpolation theorem,
\begin{equation}\label{2.1a}
\mbox{curl}\, u^\varepsilon \qquad \mbox{is compact in} \,\,\, H^{-1}_{loc}.
\end{equation}
For a fixed compact set $K$, the density $\rho^\varepsilon$ can be controlled
by $\sup\limits_{\varepsilon}\Big(\sup\limits_{x\in K} B^\varepsilon/\inf\limits_{x\in K}S^\varepsilon\Big)^{\frac{1}{\gamma-1}}$.
Similarly, we have
\begin{equation}\label{2.1b}
\mbox{div}(\rho^\varepsilon u^\varepsilon) \qquad \mbox{is compact in} \,\,\,  H^{-1}_{loc}.
\end{equation}

\medskip
{\it Step 3. The strong convergence of $u^\varepsilon(x)$, which also leads to
the strong convergence of $(\rho^\v, p^\v)(x)$ from Step 1}.
By the div-curl lemma of Murat \cite{Murat} and Tartar \cite{Tartar}, the Young measure representation theorem for a uniformly bounded sequence of functions ({\it cf.} Tartar \cite{Tartar}; also
Ball \cite{Ball}),
and \eqref{2.1a}--\eqref{2.1b},
we have the following commutation identity:
\begin{eqnarray}\label{3.3}
\sum\limits_{i=1}^n\langle \nu( u), \,\rho(q) u_i\rangle \langle \nu( u), u_i\rangle&=&\langle \nu(u),
\,\sum\limits_{i=1}^n\rho(q) u_i u_i\rangle.
\end{eqnarray}
Here and hereafter, for simplicity of notation, we have used that $\nu(u):=\nu_{x}(\rho(q; \overline{B}(x), \overline{S}(x)), u)$
denotes the associated Young measure (a probability measure)
for the sequence  $u^\varepsilon(x)$,
and $\rho(q):=\rho(q; \overline{B}(x), \overline{S}(x))$ for the limits $\overline{B}(x)$ and
$\overline{S}(x)$ of $B^\varepsilon(x)$ and $S^\varepsilon(x)$ respectively.

Then the main point in this step for the
compensated compactness framework is to prove that
$\nu$ is in fact a Dirac measure, which
in turn implies the compactness of the sequence
$u^\varepsilon(x)$.

Combining both sides of $(\ref{3.3})$ together, we have
\begin{equation}
\langle \nu(u^{(1)})\otimes\nu(u^{(2)}), \,\sum\limits_{i=1}^n\rho(q^{(1)}) u^{(1)}_i (u_i^{(1)}-u_i^{(2)})\rangle=0.
\end{equation}
Exchanging indices $(1)$ and $(2)$, we obtain the following symmetric commutation identity:
\begin{equation}\label{3.4}
\langle \nu(u^{(1)})\otimes\nu(u^{(2)}),
I(u^{(1)}, u^{(2)})\rangle=0,
\end{equation}
where
\begin{equation}\label{3.5}
I(u^{(1)},u^{(2)})
=\sum_{i=1}^n(\rho(q^{(1)}) u_i^{(1)}-\rho(q^{(2)})u_i^{(2)})(u_i^{(1)}-u_i^{(2)}).
\end{equation}

\noindent
Then it remains to prove the strong convergence of the velocity $u^\v$ from the above identity.

Notice that
\begin{eqnarray}\label{3.6}
 I(u^{(1)}, u^{(2)})&=&\sum_{i=1}^n\big(\rho(q^{(1)}) u_i^{(1)}-\rho(q^{(2)})u_i^{(2)}\big)\big(u_i^{(1)}-u_i^{(2)}\big)\nonumber\\
&=&\sum_{i=1}^n \big(\rho(q^{(1)}) |u_i^{(1)}|^2-\rho(q^{(1)}) u_i^{(1)}u_i^{(2)}-\rho(q^{(2)}) u_i^{(2)}u_i^{(1)}-\rho(q^{(2)}) |u_i^{(2)}|^2\big) \nonumber\\
&=&\rho(q^{(1)})\big((q^{(1)})^2-\sum\limits_{i=1}^n
u_i^{(1)}u_i^{(2)}\big)+\rho(q^{(2)})\big((q^{(2)})^2-\sum\limits_{i=1}^n
u_i^{(1)}u_i^{(2)}\big).
\end{eqnarray}
The Cauchy inequality implies
\begin{eqnarray}\label{3.7}
I(u^{(1)}, u^{(2)})&\geq &\rho(q^{(1)})\big((q^{(1)})^2-q^{(1)} q^{(2)}\big)+\rho(q^{(2)})\big((q^{(2)})^2-q^{(1)}q^{(2)}\big)\nonumber\\
&=&\big(q^{(1)}-q^{(2)}\big)\big(\rho(q^{(1)})q^{(1)}-\rho(q^{(2)})q^{(2)}\big)\nonumber\\
&=&\big(q^{(1)}-q^{(2)}\big)^2 \frac{d(\rho q)}{dq}(\tilde{q}),
\end{eqnarray}
where $\tilde{q}$ lies between $q^{(1)}$ and $q^{(2)}$ by the
mean-value theorem.

Taking derivative with respect to $q$ on $(\ref{3.3f})$, we obtain
$$
\frac{d \rho}{d q}=-\frac{q}{(\gamma-1)\rho^{\gamma-2}\overline{S}}=-\frac{\rho q}{c^2}.
$$
Then
$$
\frac{d(\rho q)}{d q}=\rho(1-M^2).
$$
For subsonic-sonic flows, {\it i.e.}, $q^{(1)}, q^{(2)}\leq q_{\rm cr}(\overline{B})$, we have
$$
M^2(\tilde{q})\le 1.
$$
On the other hand,
$$
\rho(\tilde{q})\geq \left(\frac{2\overline{B}}{(\gamma+1)\overline{S}}\right)^{\frac{1}{\gamma-1}}\geq 0.
$$
Notice that $M^2(\tilde{q})=1$ if and only if $q^{(1)}=q^{(2)}=q_{\rm cr}(\overline{B})$,
while $\rho(\tilde{q})=0$ if and only if $\overline{B}=0$ and $q^{(1)}=q^{(2)}=q_{\rm cr}(\overline{B})$.
Then
\begin{equation}\label{3.8}
I(u^{(1)}, u^{(2)})\ge\big(q^{(1)}-q^{(2)}\big)^2 \, \rho(\tilde{q})\big(1-M^2(\tilde{q})\big)\geq 0.
\end{equation}
With $(\ref{3.4})$, it implies that
$$
q^{(1)}=q^{(2)},
$$
which also deduces
$$
\rho^{(1)}=\rho^{(2)}, \qquad p^{(1)}=p^{(2)}.
$$

Again, using $(\ref{3.4})$ and $(\ref{3.5})$, we further obtain
\begin{eqnarray}\label{3.10}
0&=&\langle \nu(u^{(1)})\otimes\nu(u^{(2)}), I(u^{(1)}, u^{(2)})\rangle \nonumber\\
&=&\langle \nu(u^{(1)})\otimes\nu(u^{(2)}), \rho(q^{(1)})\sum\limits_{i=1}^n\big(u_i^{(1)}-u_i^{(2)}\big)^2\rangle,
\end{eqnarray}
which immediately implies
$$
u^{(1)}=u^{(2)},
$$
{\it i.e.}, $\nu(u)$ concentrates on a single point.
If this would not be the case, we could suppose to have two different points $\acute{u}$ and $\grave{u}$ in the support of $\nu$.
Then $(\acute{u}, \acute{u})$, $(\acute{u}, \grave{u})$, $(\grave{u}, \acute{u})$, and $(\grave{u}, \grave{u})$
would be in the support of $\nu\otimes\nu$, which contradicts with $u^{(1)}=u^{(2)}$.
Therefore, the Young measure $\nu$ is a Dirac measure, which implies the strong convergence
of $(\rho^\v, u^\v, p^\v)$. This completes the proof.

\medskip
For the homentropic case, the entropy function $S$ is constant.
Then the pressure $p$ may be regarded as a  function
of the density $\rho$ in
$C^1\left(\mathbb{R}^+\cup\{0\}\right)\cap C^2\left(\mathbb{R}^+\right)$,
which satisfies
$p'(\rho)>0$ and $2p'(\rho) + \rho p''(\rho) > 0$ for $\rho>0$.

Without loss of generality, we set the enthalpy $h(\rho)$ as
$$
h(\rho):=\int_{1}^{\rho}\frac{p'(s)}{s}ds,
$$
so that $h(1)=0$.
It is noticeable that, in this case, we do not have the property that the Bernoulli function is greater or equal to zero.
To replace the nonnegative property, we introduce a lower bound
$$
B_{\rm min}=\lim_{\rho\to 0^+}f(\rho),
$$
while
$$
f(\rho):=\frac{p'(\rho)}{2}+h(\rho)
$$
belongs to $C^1(\mathbb{R}^{+}\cup\{0\})$.
Since $2p'(\rho) + \rho p''(\rho)> 0$ for $\rho>0$, $f(\rho)$ is an increase function in $\rho>0$.
Then $f(\rho)>B_{\rm min}$ for $\rho>0$ and $B_{\rm min}$ is lower bound so that $B_{\rm min}<f(1)=\frac{p'(1)}{2}$.

In the homentropic case, conditions (A.1)--(A.3) can be reformulated as:

\medskip
(B.1). $M^\varepsilon\leq 1$ {\it a.e.} in
$\Omega$;

\vspace{2mm}
(B.2). $B^\varepsilon= \frac{(q^\varepsilon)^2}{2}+\int_{1}^{\rho^\varepsilon}\frac{p'(s)}{s}ds$
are uniformly bounded and $B^\varepsilon \ge B_{\rm min}$;

\vspace{2mm}
(B.3). $\mbox{curl}\  u^\varepsilon$, $\frac{e_1(\varepsilon)}{\rho^\varepsilon}$, and $\frac{e_2(\varepsilon)}{\rho^\varepsilon}$ are
uniformly bounded measures.

\medskip
Similar to Theorem $\ref{thm3.1}$, we have

\begin{theorem}[Compensated  compactness framework for the homentropic case]\label{thm2.1}
Let a sequence of functions $\rho^\varepsilon(x)$ and $u^\varepsilon(x)=(u^\varepsilon_1,
\cdots, u^\varepsilon_n)(x)$ satisfy conditions {\rm (B.1)--(B.3)}.
Then there exists a subsequence (still labeled)  $(\rho^\varepsilon, u^\varepsilon)(x)$
such that
$$
\rho^\varepsilon(x)\rightarrow \rho(x),
\quad
u^\varepsilon(x)\rightarrow u(x)=(u_1, \cdots, u_n)(x)
 \qquad \mbox{a.e. $x\in\Omega\, \,$
as $\varepsilon\rightarrow 0$}
$$
and
$$
M(x)\leq 1 \qquad \mbox{a.e. $x\in \Omega$.}
$$
\end{theorem}

\noindent\textbf{Proof}.
First, since $p'(\rho)>0$ for $\rho>0$, we can employ the implicit function theorem
to conclude $\rho(q; B)=h^{-1}(B-\frac{q^2}{2})$.
Then we can regard $\rho^\varepsilon$ as a function of $q^\varepsilon$ and $B^\varepsilon$,
that is, $\rho^\varepsilon=\rho(q^\varepsilon; B^\varepsilon)$.
As a consequence, the sequence $\rho^\varepsilon$ is nonnegative and uniformly bounded.
Conditions (B.1)--(B.2) indicate directly that the speed sequence $q^\varepsilon$ is uniformly bounded.

Differentiating the Bernoulli functions $B^\varepsilon$ with respect to $x_i$ yields
\begin{equation}\label{2.2}
\begin{array}{ll}
\partial_{x_i} B^\varepsilon &=\sum\limits_{j=1}^n u_j^\varepsilon\, \partial_{x_i} u_j^\varepsilon+\frac{p'(\rho^\varepsilon)}{\rho^\varepsilon}\partial_{x_i}\rho^\varepsilon\\
& = \sum\limits_{j=1}^n u_j^\varepsilon\,\partial_{x_i}\, u_j^\varepsilon
 + \frac{\partial_{x_i} p(\rho^\varepsilon)}{\rho^\varepsilon}\\
& = \sum\limits_{j=1}^n u_j^\varepsilon\, \partial_{x_i} u_j^\varepsilon + \frac{e_{2i}(\varepsilon)-\sum\limits_{j=1}^n\rho^\varepsilon u_j^\varepsilon\, \partial_{x_j} u_i^\varepsilon-u_i^\varepsilon e_1(\varepsilon)}{\rho^\varepsilon}\\
& = \sum\limits_{j=1}^n u_j^\varepsilon (\partial_{x_i} u_j^\varepsilon-\partial_{x_j} u_i^\varepsilon)+ \frac{e_{2i}(\varepsilon)-u_i^\varepsilon e_1(\varepsilon)}{\rho^\varepsilon}\\
& = \sum\limits_{j=1}^n u_j^\varepsilon\, \omega^\varepsilon_{ij}+ \frac{e_{2i}(\varepsilon)-u_i^\varepsilon e_1(\varepsilon)}{\rho^\varepsilon}.
\end{array}
\end{equation}

From the boundedness of $q^\varepsilon$ and condition (B.3), we conclude that
$$
\nabla B^\varepsilon \qquad\mbox{is uniformly bounded measures}.
$$
Since $B^\varepsilon$ is uniformly bounded, the total-variation norm of $B^\varepsilon$
is uniformly bounded, which implies that
$B^\varepsilon(x)$ converges to $\overline{B}(x)$ in $L^1_{loc}$, as $\varepsilon$ tends to $0$,
and $\overline{B}(x)\ge B_{\rm min}$ {\it a.e.} in $\Omega$.

\smallskip
From the definition,
$(\mbox{curl}\, u^\varepsilon)_{ij}=\omega_{ij}^\varepsilon=\partial_{x_i} u_j^\varepsilon - \partial_{x_j} u_i^\varepsilon$ is bounded
in $W^{-1, \infty}$. On the other hand, $\mbox{curl}\, u^\varepsilon$ is a uniformly bounded measure sequence, which implies
that
$\mbox{curl}\, u^\varepsilon$ is compact in $W^{-1, q}$ for each $1\le p <\frac{n}{n-1}$.
From the interpolation theorem,
$$
\mbox{curl}\, u^\varepsilon \qquad \mbox{is compact in}~ H^{-1}_{loc}.
$$
Similarly, we have
$$
\mbox{div}(\rho^\varepsilon u^\varepsilon) \qquad \mbox{is compact in}~ H^{-1}_{loc}.
$$

Next, we will discuss the strong convergence in two cases.

\medskip
{\it Case 1: $\overline{B}>B_{\rm min}$}.
From the monotonicity of $f(\rho)$,
we have
$$
\rho_{\rm cr}>0,
$$
which implies $f(\rho_{\rm cr})=\overline{B}$ and, on the support of $\nu$,
$$
\rho(q)\in[\rho_{\rm cr}, h^{-1}(\overline{B})].
$$
Following similar argument for Theorem $\ref{thm3.1}$, we obtain the following commutation identity:
\begin{equation}\label{2.41}
\langle \nu(u^{(1)})\otimes\nu(u^{(2)}), \,
I( u^{(1)}, u^{(2)})\rangle=0,
\end{equation}
and
\begin{equation}\label{2.51}
I(u^{(1)}, u^{(2)})
=\sum_{i=1}^n\big(\rho(q^{(1)}) u_i^{(1)}-\rho(q^{(2)})u_i^{(2)}\big)\big(u_i^{(1)}-u_i^{(2)}\big).
\end{equation}

By the same calculation as in the previous argument, we have
\begin{equation}\label{2.7}
I(u^{(1)}, u^{(2)})\geq\big(q^{(1)}-q^{(2)}\big)^2 \frac{d(\rho q)}{dq}(\tilde{q}),
\end{equation}
where $\tilde{q}$ lies between $q^{(1)}$ and $q^{(2)}$.

From the definition of $\rho(q)$, we obtain that $\frac{d \rho}{d q}=-\frac{\rho q}{c^2}$,
which implies $\frac{d(\rho q)}{dq}=\rho(1-M^2)$.
From subsonic-sonic flows, {\it i.e.}, $q^{(1)}, q^{(2)}\leq q_{\rm cr}$,
we have
$$
\rho(\tilde{q})\big(1-M^2(\tilde{q})\big)\geq 0.
$$
Notice that $M^2(\tilde{q})=1$
if and only if $q^{(1)}=q^{(2)}=q_{\rm cr}$.
Then
\begin{equation}\label{2.8}
I(u^{(1)}, u^{(2)})\ge\left(q^{(1)}-q^{(2)}\right)^2 \, \rho(\tilde{q})\big(1-M^2(\tilde{q})\big)\geq 0.
\end{equation}
Thus, from $(\ref{2.41})$, we obtain
$$
q^{(1)}=q^{(2)}.
$$

Again using $(\ref{2.41})$ and $(\ref{2.51})$, we have
\begin{eqnarray}\label{3.10a}
0&=&\langle \nu(u^{(1)})\otimes\nu(u^{(2)}), \, I(u^{(1)}, u^{(2)})\rangle \nonumber\\
&=&\langle \nu(u^{(1)})\otimes\nu(u^{(2)}), \, \rho(q^{(1)})\sum\limits_{i=1}^n\big(u_i^{(1)}-u_i^{(2)}\big)^2\rangle,
\end{eqnarray}
which immediately implies $u^{(1)}=u^{(2)}$, {\it i.e.},
the Young measure $\nu$ is a Dirac measure.
Thus, we conclude the strong convergence
of $(\rho^\varepsilon, u^\varepsilon)(x)$ to $(\rho, u)(x)$
with $M(x)\le 1$ a.e $x\in \Omega$.

\medskip
{\it Case 2. $\overline{B}=B_{\rm min}$}.
Considering the boundedness of $\overline{B}$, we can see $\lim\limits_{\rho\rightarrow0^+}\int_{1}^{\rho}\frac{p'(s)}{s}ds$
is finite, which implies $p'(0)=0$. On the other hand, the monotonicity of $f(\rho)$ shows
that  $\rho(q)=0$ for any $q$ belonging to the support of the Young measure $\nu$.
It follows from the subsonic-sonic condition that $0\le q^2\le p'(0)=0$.
Then the strong convergence of $(\rho^\varepsilon, u^\varepsilon)$
is shown in this case.

\medskip
\begin{remark}\label{rem3}
The main difference between Theorems $\ref{thm3.1}$ and $\ref{thm2.1}$ is
that the homentropic case has the Bernoulli-vortex relation $(\ref{2.2})$,
while the full Euler case has
\begin{equation}\label{2.17}
\partial_{x_i} B^\varepsilon= \sum\limits_{j=1}^n u^\varepsilon_j\, \omega^\varepsilon_{ij}
+\frac{(\rho^\varepsilon)^{\gamma-1}}{\gamma}\partial_{x_i} S^\varepsilon
+\frac{e_{2i}(\varepsilon)-u_i^\varepsilon e_1(\varepsilon)}{\rho^\varepsilon}
\end{equation}
so that the gradient of $B^\varepsilon$ cannot be achieved only
by the vorticity $\omega^\varepsilon$.
\end{remark}

\begin{remark}  The main theorem in Huang-Wang-Wang \cite{Huang-Wang-Wang} is included in Theorem $\ref{thm2.1}$.
Condition {\rm (B.1)} is the same as the one in \cite{Huang-Wang-Wang},
while  $B^\varepsilon$ is assumed to be constant in \cite{Huang-Wang-Wang}
which clearly satisfies condition {\rm (B.2)}.
As a consequence, condition {\rm (B.3)} could be regarded as the $H^{-1}_{loc}$--compactness of $\mbox{curl}\, u^\varepsilon$.
The irrotational condition in \cite{Huang-Wang-Wang} is removed. Thus, Theorem $\ref{thm2.1}$ includes
more physical consideration.
From $(\ref{2.2})$, the irrotational condition implies that the Bernoulli function is a constant in the flow field.
\end{remark}

\begin{remark}\label{rem4}
Consider any function
$Q(\rho, u, p)=(Q_1,\cdots,Q_n)(\rho, u, p)$ satisfying
\begin{equation}\label{3.12}
\mbox{\rm div}\,(Q(\rho^\varepsilon, u^\varepsilon,  p^\varepsilon))=o_Q(\varepsilon),
\end{equation}
where $o_Q(\varepsilon)\rightarrow 0$ in the distributional sense as
$\varepsilon\rightarrow 0$.
We can see from the strong convergence of
$(\rho^\varepsilon, u^\varepsilon, p^\varepsilon)$ ensured by Theorem $\ref{thm3.1}$ that
$\mbox{\rm div} (Q(\rho, u, p))=0$ holds in the distributional sense.
Thus, if
\begin{equation}
\mbox{\rm div} \, (\rho^\varepsilon u^\varepsilon\otimes
u^\varepsilon+p^\varepsilon I)=e_2(\varepsilon)\rightarrow 0 \qquad
\mbox{in the sense of distributions},
\end{equation}
the weak solution also satisfies the
momentum equations in $(\ref{1.5})_2$ and the energy equation $(\ref{1.5})_3$ in the distributional sense.
The statement is also valid for Theorem $\ref{thm2.1}$.
\end{remark}

Then, as corollaries, we conclude the following theorems.

\begin{theorem}[Convergence of approximate solutions for the full Euler flow]\label{thm3.2}
Let $\rho^\varepsilon(x)$, $u^\varepsilon(x)=(u^\varepsilon_1, \cdots, u^\varepsilon_n)(x)$, and $p^\varepsilon(x)$
be a sequence of approximate solutions satisfying  {\rm (A.1)}--{\rm (A.3)} and $e_j(\varepsilon)\rightarrow 0, j=1,2,3$,
in the distributional sense as $\varepsilon\rightarrow 0$.
Then there exists a
subsequence (still labeled) $(\rho^\varepsilon, u^\varepsilon, p^\varepsilon)(x)$ that converges a.e.
as $\varepsilon\rightarrow 0$ to a weak solution $(\rho, u, p)$ to the Euler equations of $(\ref{1.5})$,
which satisfies  $M(x) \le 1$, {\it a.e.} $x\in \Omega$.
\end{theorem}

\begin{theorem}[Convergence of approximate solutions for the homentropic Euler flow]\label{thm2.2}
Let $\rho^\varepsilon(x)$ and $u^\varepsilon(x)=(u^\varepsilon_1, \cdots, u^\varepsilon_n)(x)$
be a sequence of approximate solutions satisfying  {\rm (B.1)}--{\rm (B.3)}
and $e_j(\varepsilon)\rightarrow 0, j=1,2$, in the distributional sense
as $\varepsilon\rightarrow 0$. Then there exists a
subsequence (still labeled) $(\rho^\varepsilon, u^\varepsilon)(x)$ that converges {\it a.e.}
as $\varepsilon\rightarrow 0$ to a weak solution $(\rho, u)$
to the Euler equations of $(\ref{1.1})$, which satisfies
$M(x) \le 1$ {\it a.e.} $x\in \Omega$.
\end{theorem}

There are various ways to construct approximate solutions by either
numerical methods or analytical methods such as vanishing viscosity
methods. As direct applications, we show two examples
in Sections 3--4 to apply the
compactness framework built above in establishing
existence theorems for multidimensional subsonic-sonic full Euler flows through infinitely long nozzles.

\section{Subsonic-Sonic Limit for Two-Dimensional Steady Full Euler Flows \\
in an Infinitely Long Nozzle}

In this section, as a direct application of the compactness framework
established in Theorem \ref{thm3.1}, we obtain the subsonic-sonic limit of
steady subsonic full Euler flows in a two-dimensional, infinitely long nozzle.

The infinitely long nozzle is defined as
\begin{equation*}
\Omega=\{(x_1,x_2)\, :\, f_1(x_1)<x_2<f_2(x_1), \,-\infty<x_1<\infty\},
\end{equation*}
with the nozzle walls $\partial\Omega:=W_1\cup W_2$, where
\begin{equation*}
W_i=\{(x_1,x_2)\, :\, x_2=f_i(x_1)\in C^{2,\alpha}, ~-\infty<x_1<\infty\}, \qquad i=1,2.
\end{equation*}

\begin{figure}[htbp]
\small \centering
\includegraphics[width=8.5cm]{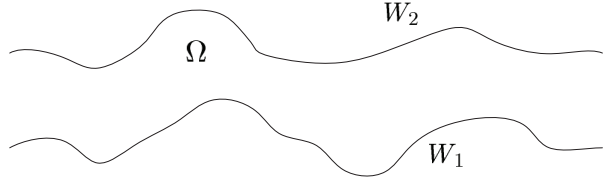}
\caption{Two-Dimensional Infinitely Long Nozzle}
\label{Fig1}
\end{figure}

Suppose that $W_1$ and $W_2$ satisfy
\begin{align}\label{cdx-3}
&f_2(x_1)>f_1(x_1) ~ \qquad \mbox{for} ~x_1\in(-\infty, \infty),\nonumber\\
&f_1(x_1)\rightarrow 0, \quad f_2(x_1)\rightarrow 1 \qquad \mbox{as} ~x_1\rightarrow -\infty, \nonumber\\
&f_1(x_1)\rightarrow a, \quad f_2(x_1)\rightarrow b>a
  \qquad \mbox{as} ~x_1\rightarrow \infty,
\end{align}
and there exists $\alpha>0$ such that
\begin{equation}\label{cdx-4}
\|f_i\|_{C^{2,\alpha}(\mathbb{R})}\leq C,
\qquad i=1,2,
\end{equation}
for some positive constant $C$.
It follows that $\Omega$ satisfies the uniform exterior sphere condition
with some uniform radius $r>0$. See Fig \ref{Fig1}.

\bigskip

Suppose that the nozzle has impermeable solid walls so that the flow satisfies
the slip boundary condition:
\begin{equation}\label{cdx-5}
u \cdot\nu =0 \qquad \mbox{on} ~\partial \Omega,
\end{equation}
where $u=(u_1, u_2)\in \R^2$ is the velocity and $\nu=(\nu_1, \nu_2)$ is
the unit outward normal to the nozzle wall.
In the flow without vacuum, it can be written as
\begin{equation}\label{cdx-sc}
(\rho u) \cdot\nu =0 \qquad \mbox{on} ~\partial \Omega.
\end{equation}

It follows from $(\ref{1.5})_1$ and $(\ref{cdx-sc})$ that
\begin{equation}\label{cdx-6}
\int_s\, (\rho u) \cdot l \, ds\equiv m
\end{equation}
holds for some constant $m$, which is the mass flux,
where $s$ is any curve transversal to the $x_1$--direction,
and $l$ is the normal of $s$ in the positive $x_1$--axis
direction.

We assume that the upstream entropy function is given, {\it i.e.},
\begin{equation}\label{cdx-8}
\frac{\gamma p}{(\gamma-1)\rho^\gamma}\longrightarrow S_-(x_2)
\qquad \mbox{as} ~x_1\rightarrow -\infty,
\end{equation}
and the upstream Bernoulli function is given, {\it i.e.},
\begin{equation}\label{cdx-11}
\frac{q^2}{2}+\frac{\gamma p}{(\gamma-1)\rho}\longrightarrow B_-(x_2)
\qquad\mbox{as} ~x_1\rightarrow -\infty,
\end{equation}
where $B(x_2)$ is a function defined on $[0,1]$.

\medskip
$\mathbf{Problem~1}(m)$:
Solve the full Euler system $(\ref{1.5})$ with the boundary condition $(\ref{cdx-sc})$,
the mass flux condition $(\ref{cdx-6})$,
and the asymptotic conditions $(\ref{cdx-8})$--$(\ref{cdx-11})$.

Set
\begin{equation*}
\underline{S}=\inf_{x_2\in[0,1]}S_-(x_2),
\qquad \underline{B}=\inf_{x_2\in[0,1]}B_-(x_2).
\end{equation*}

For this problem, the following theorem has been established in
Chen-Deng-Xiang \cite{ChenDX}.

\begin{theorem}\label{thm5.1} Let the nozzle walls $\partial \Omega$ satisfy
$(\ref{cdx-3})$--$(\ref{cdx-4})$, and let
$\underline{S}>0$ and $\underline{B}>0$. Then there exists
$\delta_0>0$ such that, if $\|(S_--\underline{S},B_--\underline{B})\|_{C^{1,1}([0,1])}\leq \delta$
for $0<\delta\leq\delta_0$, $(S_-B_-^{-\gamma})'(0)\geq 0$, and $(S_-B_-^{-\gamma})'(1)\leq 0$,
there exists $\hat{m}\geq 2\delta_0^{{1}/{8}}$ such that,
for any $m\in(\delta^{{1}/{4}},\hat{m})$,
there exists a global solution ({\it i.e.} a full Euler flow)
$(\rho, u, p) \in C^{1,\alpha}(\overline{\Omega})$
of $\mathbf{Problem~1}(m)$ such that the following hold:

\medskip
{\rm (i)} Subsonicity and positivity of the horizontal velocity: The flow is uniformly
subsonic
with positive horizontal velocity in the whole nozzle, {\it i.e.},
\begin{equation}\label{cdx3}
\sup_{\overline{\Omega}}(q^2-c^2)<0, \quad u_1>0  \qquad \mbox{in} ~\overline{\Omega};
\end{equation}

{\rm (ii)} The flow satisfies the following asymptotic behavior in the far fields:
As $x_1\rightarrow -\infty$,
\begin{equation}\label{cdx4}
p\rightarrow p_->0, \qquad
u_1\rightarrow u_-(x_2)>0, \qquad
(u_2,\rho)\rightarrow (0,\rho_-(x_2;p_-)),
\end{equation}
\begin{equation}\label{cdx5}
\nabla p\rightarrow 0, \qquad
\nabla u_1\rightarrow (0, u_-'(x_2)),
\qquad \nabla u_2\rightarrow 0,
\qquad \nabla\rho\rightarrow (0,\rho_-'(x_2;p_-))
\end{equation}
uniformly for $x_2\in K_1\Subset(0,1)$, where
$
\rho_-(x_2;p_-)=\big(\frac{\gamma p_-}{(\gamma-1)S_-(x_2)}\big)^{{1}/{\gamma}},
$
the constant $p_-$ and function $u_-(x_2)$ can be determined by $m$, $S_-(x_2)$, and $B_-(x_2)$ uniquely;

\smallskip
{\rm (iii)} Uniqueness: The full Euler flow of $\mathbf{Problem~1} (m)$ satisfying $(\ref{cdx3})$
and the asymptotic behavior $(\ref{cdx4})$--$(\ref{cdx5})$ is unique.

\smallskip
{\rm (iv)} Critical mass flux: $\hat{m}$ is the upper critical mass flux for the existence
of subsonic flow in the following sense: Either $\sup\limits_{\overline{\Omega}}(q^2-c^2)\rightarrow 0$
as $m\rightarrow \hat{m}$,
or there is no $\sigma>0$ such that, for all
$m\in(\hat{m},\hat{m}+\sigma)$, there are full Euler flows of $\mathbf{Problem~1} (m)$
satisfying $(\ref{cdx3})$, the asymptotic behavior $(\ref{cdx4})$--$(\ref{cdx5})$, and
$\sup\limits_{m\in(\hat{m},\hat{m}+\sigma)}
\sup\limits_{\overline{\Omega}}(c^2-q^2)>0$.
\end{theorem}

We note that Theorem \ref{thm5.1} does not apply to the critical flows,
that is, those flows for which $m=\hat{m}$  must be sonic at some point.
Now we can employ Theorem $\ref{thm3.1}$
to establish a more general result.

\begin{theorem}[Subsonic-sonic limit of two-dimensional full Euler flows]\label{thm5.2}
Let $\delta^{{1}/{4}}<m^{\varepsilon}<\hat{m}$ be a sequence of mass fluxes,
and let $(\rho^\varepsilon, u^\varepsilon, p^\varepsilon)(x)$ be the corresponding
sequence of solutions to $\mathbf{Problem~1} (m^\varepsilon)$.
Then, as $m^{\varepsilon}\rightarrow \hat{m}$, the solution
sequence possesses a subsequence (still denoted by)
$(\rho^\varepsilon, u^{\varepsilon}, p^\varepsilon)(x)$
that converges strongly {\it a.e.} in $\Omega$ to
a vector function $(\rho, u, p)(x)$ which is a weak solution
of $\mathbf{Problem~1} (\hat{m})$.
Furthermore, the limit solution $(\rho, u, p)(x)$ also satisfies $(\ref{1.5})$
in the distributional sense and the boundary conditions $(\ref{cdx-sc})$
as the normal trace of the divergence-measure field
$\rho u$ on the boundary in the sense of Chen-Frid {\rm \cite{Chen7}}.
\end{theorem}

\noindent{\bf Proof}.  We divide the proof into three steps.

\smallskip
1. We first need to show that $(\rho^\varepsilon, u^\varepsilon, p^\varepsilon)(x)$
satisfy condition (A.1)--(A.3).
For $B^\varepsilon$ and $S^\varepsilon$, we have
\begin{eqnarray}\label{3.2.1}
\begin{cases}
\partial_{x_1}(\rho^\varepsilon u_1^\varepsilon)+\partial_{x_2}(\rho^\varepsilon u_2^\varepsilon)=0,\\[2mm]
\partial_{x_1}(\rho^\varepsilon u_1^\varepsilon B^\varepsilon)+\partial_{x_2}(\rho^\varepsilon u_2^\varepsilon B^\varepsilon)=0,\\[2mm]
\partial_{x_1}(\rho^\varepsilon u_1^\varepsilon S^\varepsilon)+\partial_{x_2}(\rho^\varepsilon u_2^\varepsilon S^\varepsilon)=0.
\end{cases}
\end{eqnarray}

From $(\ref{3.2.1})_1$, we introduce the following stream function $\psi^\varepsilon$:
\begin{eqnarray}
\begin{cases}
\partial_{x_1}\psi^\varepsilon=-\rho^\varepsilon u_2^\varepsilon,\\[2mm]
\partial_{x_2}\psi^\varepsilon=\rho^\varepsilon u_1^\varepsilon,
\end{cases}
\end{eqnarray}
which means that $\psi^\varepsilon$ is constant along the streamlines.

From the far-field behavior of the Euler flows,
we define 
$$
\psi^\varepsilon_-(x_2):=\lim\limits_{x_1\rightarrow-\infty}\psi^\varepsilon(x_1, x_2).
$$
Since both the upstream Bernoulli and entropy functions are given,
$B^\varepsilon$ and $S^\varepsilon$ have the following expression:
$$
B^\varepsilon(x)=B_{-}((\psi^\varepsilon_{-})^{-1}(\psi^\varepsilon(x))),
\qquad S^\varepsilon(x)=S_{-}((\psi^\varepsilon_{-})^{-1}(\psi^\varepsilon(x))),
$$
where $(\psi^\varepsilon_{-})^{-1}\psi^\varepsilon(x)$ is a function
from $\Omega$ to $[0,1]$.
For fixed $x_1$, it can be regarded as a backward characteristic map with
$$
\frac{\partial ((\psi^\varepsilon_{-})^{-1}\psi^\varepsilon)}{\partial x_2}
=\frac{\rho^\varepsilon u_1^\varepsilon}{\rho^\varepsilon_-u^\varepsilon_-}>0.
$$
The boundedness and positivity of $\rho^\varepsilon_-u_-^\varepsilon$
and  $\rho^\varepsilon u_1^\varepsilon$ show that
the map is not degenerate.
Thus, we have
\begin{eqnarray}
\begin{cases}
\partial_{x_1}B^\varepsilon(x)=-B'_{-}((\psi^\varepsilon_{-})^{-1}(\psi^\varepsilon(x)))\frac{\rho^\varepsilon u^\varepsilon_2}{\rho^\varepsilon_-u^\varepsilon_{-}},\\[2mm]
\partial_{x_2}B^\varepsilon(x)=B'_{-}((\psi^\varepsilon_{-})^{-1}(\psi^\varepsilon(x)))\frac{\rho^\varepsilon u^\varepsilon_1}{\rho^\varepsilon_-u^\varepsilon_-}.
\end{cases}
\end{eqnarray}
Then $B^\varepsilon$ is uniformly bounded in $BV$, which implies its strong convergence.
The similar argument can lead to the strong convergence of $S^\varepsilon$.

2. For the corresponding vorticity sequence $\omega^\varepsilon$, (\ref{2.17}) can be written as
\begin{eqnarray}
\begin{cases}
\partial_{x_1} B^\varepsilon=  u_2^\varepsilon \omega^\varepsilon+\frac{1}{\gamma}(\rho^\varepsilon)^{\gamma-1}\partial_{x_1} S^\varepsilon,\\[2mm]
\partial_{x_2} B^\varepsilon=  - u_1^\varepsilon \omega^\varepsilon+\frac{1}{\gamma}(\rho^\varepsilon)^{\gamma-1}\partial_{x_2} S^\varepsilon.
\end{cases}
\end{eqnarray}
By direct calculation, we have
\begin{eqnarray}
\omega^\varepsilon&=&\frac{1}{(q^\varepsilon)^2}
\Big(u_2^\varepsilon\big(\partial_{x_1} B^\varepsilon-\frac{1}{\gamma}(\rho^\varepsilon)^{\gamma-1} \partial_{x_1} S^\varepsilon\big)
-u_1^\varepsilon\big(\partial_{x_2} B^\varepsilon-\frac{1}{\gamma}(\rho^\varepsilon)^{\gamma-1}\partial_{x_2} S^\varepsilon\big)\Big)\nonumber\\
&=&\frac{1}{\rho_-^\varepsilon u_-^\varepsilon}\Big(\rho^\varepsilon B_-'-\frac{1}{\gamma}(\rho^\varepsilon)^\gamma S_-'\Big),
\end{eqnarray}
which implies that $\omega^\varepsilon$ as a measure sequence is uniformly bounded, which is compact in $H^{-1}_{loc}$.

\smallskip
Then Theorem $\ref{thm3.1}$ immediately implies that
the solution sequence has a subsequence (still denoted
by) $(\rho^\varepsilon, u^{\varepsilon}, p^\varepsilon)(x)$
that converges {\it a.e.} in $\Omega$ to
a vector function $(\rho, u, p)(x)$.

\smallskip
Since
$(\ref{1.5})$ holds for the sequence of subsonic solutions
$(\rho^\varepsilon, u^{\varepsilon}, p^\varepsilon)(x)$,
it is straightforward to see that $(\rho, u, p)$ also
satisfies $(\ref{1.5})$ in the distributional sense.

\smallskip
3. The boundary condition is satisfied
in the sense of Chen-Frid \cite{Chen7}, which  implies
\begin{equation}
\int_{\partial \Omega}\phi(w)(\rho u)(w)\cdot \nu (w)\, d\mathcal{H}^{1}(w)
= \int_{\Omega} (\rho u)(x) \cdot \nabla \phi(x)\, dx + \langle \mbox{div}(\rho u)|_{\Omega}, \phi\rangle
\end{equation}
for $\psi\in\mathbf{C}^1_0$.
From above, we can see that $\langle \mbox{div}(\rho u)|_{\Omega}, \phi\rangle=0$. Also,
\begin{equation}
\int_{\Omega}(\rho u)(x) \cdot \nabla \phi(x)\, dx
=\lim\limits_{\varepsilon \rightarrow 0}\int_{\Omega}(\rho^\varepsilon u^\varepsilon)(x) \cdot \nabla \phi(x)\, dx=0.
\end{equation}
Then we have
\begin{equation}
\int_{\partial \Omega}\phi(w)(\rho u)(w)\cdot\nu(w)\, d\mathcal{H}^{1}(w)=0,
\end{equation}
that is, $(\rho u)\cdot \nu = 0$ on $\partial \Omega$ in $\mathcal{D}'$.

This completes the proof.

\section{Subsonic-Sonic Limit for the Full Euler Flows \\ in an Infinitely Long Axisymmetric Nozzle}

\begin{figure}[htbp]
\small \centering
\includegraphics[width=5cm]{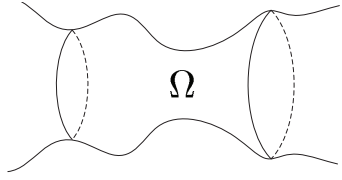}
\caption{Infinitely Long Axisymmetric Nozzle}
 \label{Fig3}
\end{figure}

We consider flows though an infinitely long axisymmetric nozzle given by
\begin{equation*}
\Omega=\{(x_1,x_2,x_3)\in\mathbb{R}^3\, : \, 0\leq \sqrt{x_2^2+x_3^2}<f(x_1),~
-\infty<x_1< \infty\},
\end{equation*}
where $f(x_1)$ satisfies
\begin{align}\label{dd-2}
&f(x_1)\rightarrow 1 \qquad \mbox{as} ~x_1\rightarrow-\infty,\nonumber\\
&f(x_1)\rightarrow r_0  \qquad \mbox{as} ~x_1\rightarrow \infty,\nonumber\\
&\|f\|_{C^{2,\alpha}(\mathbb{R})}\leq C \qquad \mbox{for}
~\mbox{some} ~\alpha>0 ~ \mbox{and} ~C>0, \\
& \inf_{x_1\in\mathbb{R}}f(x_1)=b>0.
\end{align}
See  Fig. \ref{Fig3}.

\bigskip

\medskip
The boundary condition is set as follows:
Since the nozzle wall is solid, the flow satisfies the slip boundary condition:
\begin{equation}\label{dL-3}
u\cdot\nu=0 \qquad \mbox{on} ~\partial \Omega,
\end{equation}
where $u=(u_1,u_2,u_3)$, and $\nu=(\nu_1, \nu_2, \nu_3)$ is the unit outward
normal to the nozzle wall.
In the flow without vacuum, it can be written as
\begin{equation}\label{dL-sc}
(\rho u) \cdot\nu =0  \qquad \mbox{on} ~\partial \Omega.
\end{equation}
The continuity equation in $(\ref{1.5})_1$ and the
boundary condition $(\ref{dL-sc})$ imply that the mass flux
\begin{equation}\label{dd-4}
\int_\Sigma \, (\rho u)\cdot l\, ds\equiv m_0
\end{equation}
remains for some positive constant $m_0$,
where $\Sigma$ is any surface transversal to the $x_1$--axis direction,
and $l$ is the normal of $\Sigma$ in the positive $x_1$--axis direction.

\smallskip
In Duan-Luo \cite{Duan-Luo},
the axisymmetric flows without swirl are considered for
the fluid density $\rho=\rho(x_1,r)$,
the velocity
$$
u=(u_1, u_2, u_3)=(U(x_1,r), V(x_1, r)\frac{x_2}{r}, V(x_1, r)\frac{x_3}{r}),
$$
and the pressure $p=p(x_1, r)$
in the cylindrical coordinates,
where $u_1$, $u_2$, $u_3$ are the axial velocity, radial velocity,
and swirl velocity, respectively,
and $r=\sqrt{x_2^2+x_3^2}$.
Then, instead of $(\ref{1.5})$, we have
\begin{equation}\label{dd-6}\left\{ \begin{split}
&\partial_{x_1}(r\rho U)+\partial_r(r\rho V)=0,\\
&\partial_{x_1}(r\rho U^2)+\partial_r(r\rho UV)_r+r\partial_{x_1}p=0,\\
&\partial_{x_1}(r\rho UV)+\partial_r(r\rho V^2)_r+r\partial_rp=0,\\
&\partial_{x_1}(r \rho U (E+\frac{p}{\rho}))+\partial_r(r \rho V (E+\frac{p}{\rho}))=0.
\end{split} \right.\end{equation}
Rewrite the axisymmetric nozzle as
\begin{equation*}
\Omega=\{(x_1,r)\, :\, 0\leq r<f(x_1), ~-\infty<x_1<\infty\}
\end{equation*}
with the
boundary of the nozzle:
\begin{align*}
\partial \Omega=\{(x_1,r)\, :\, r=f(x), ~-\infty<x_1<\infty\}.
\end{align*}

\smallskip
The boundary condition $(\ref{dL-3})$ becomes
\begin{equation}\label{dd-7}
(U,V,0)\cdot \tilde{\nu}=0 \qquad \mbox{on} ~\partial \Omega,
\end{equation}
where $\tilde{\nu}$ is the unit outer normal of the nozzle walls
in the cylindrical coordinates.
The mass flux condition $(\ref{dd-4})$ can
be rewritten in the cylindrical coordinates as
\begin{equation}\label{dd-9}
\int_\Sigma (r\rho U, r \rho V,0)\cdot {l}\, dS\equiv m:=\frac{m_0}{2\pi},
\end{equation}
where $\Sigma$ is any curve transversal to the $x_1$-axis direction,
and ${l}$ is the unit normal of $\Sigma$.

\smallskip
The quantities
$B=h(\rho, p)+\frac{U^2+V^2}{2}$
and $S=\frac{\gamma p}{(\gamma-1) \rho^\gamma}$ are both constants along each streamline.
For the full Euler flows in the axisymmetric nozzle,
we assume that the upstream Bernoulli and entropy functions are given, that is,
\begin{equation}\label{dd-13}
h(\rho,p)+\frac{U^2+V^2}{2}\longrightarrow B_-(r) \qquad \mbox{as} ~x_1\rightarrow -\infty,
\end{equation}
\begin{equation}\label{dd-13s}
\frac{\gamma p}{(\gamma-1) \rho^\gamma}\longrightarrow S_-(r) \qquad \mbox{as} ~x_1\rightarrow -\infty,
\end{equation}
where $B_-(r)$ and $S_-(r)$ are smooth functions defined on $[0, 1]$.

\smallskip
Set
\begin{equation}
\underline{B}=\inf\limits_{r\in[0,1]} B_-(r),
\quad \sigma_1=\|B'_-\|_{C^{0,1}([0,1])},
\end{equation}
\begin{equation}
\underline{S}=\inf\limits_{r\in[0,1]} S_-(r),
\quad \sigma_2=\|S'_-\|_{C^{0,1}([0,1])}.
\end{equation}
We denote the above problem as $\mathbf{Problem~2} (m)$. It is shown in \cite{Duan-Luo} that

\begin{theorem}\label{thm4.5}
Suppose that the nozzle satisfies $(\ref{dd-2})$.
Let the upstream Bernoulli function $B_-(r)$ and entropy function $S_-(r)$
satisfy $\underline{B}>0$, $B_{-}'(r)\in C^{1,1}([0,1])$, $B_{-}'(0)=0$, $B_{-}'(r)\geq 0$ on $r\in [0,1]$;
and $\underline{S}>0$, $S_{-}'(r)\in C^{1,1}([0,1])$, $S_{-}'(0)=0$, $S_{-}'(r)\leq 0$ on $r\in [0,1]$.
Then 

\begin{enumerate}
\item[(i)]
There exists $\delta_0>0$ such that, if $\delta:=\max\{\sigma_1, \sigma_2\}\leq\delta_0$,
then there is $\hat{m}\leq 2\delta_0^{{1}/{8}}$ so that,
for any $m\in (\delta^{{1}/{4}},\hat{m})$,
there exists a global $C^1$--solution ({\it i.e.} a full Euler flow)
$(\rho, U, V, p)\in C^1(\overline{\Omega})$ through the nozzle with mass flux
condition $(\ref{dd-9})$ and the upstream asymptotic condition $(\ref{dd-13})$.
Moreover, the flow is uniformly subsonic, and the axial velocity is always
positive, {\it i.e.},
\begin{equation}\label{dd2}
\sup_{\overline{\Omega}}(U^2+V^2-c^2)<0 \quad \mbox{and} \quad U>0
\qquad \mbox{in} ~\overline{\Omega}.
\end{equation}

\item[(ii)] The subsonic flow satisfies the following properties: As $x_1\rightarrow-\infty$,
\begin{align}\label{dd4}
&\rho\rightarrow\rho_->0, \qquad \nabla\rho\rightarrow 0, \qquad p\rightarrow \frac{\gamma-1}{\gamma}S_-(r)\rho_-^\gamma,
\qquad \nabla p\rightarrow (0, \frac{\gamma-1}{\gamma}S_-'(r)\rho_-^\gamma), \nonumber\\
&(U, V)\rightarrow (U_-(r), 0), \qquad \nabla U\rightarrow(0,U_-'(r)), \qquad \nabla V\rightarrow 0
\end{align}
uniformly for $r\in K_1\Subset(0,1)$,
where $\rho_-$ is a positive constant, and $\rho_-$ and $U_-(r)$ can
be determined by $m$, $B_-(r)$, and $S_-(r)$ uniquely.

\item[(iii)]
There exists at most one smooth axisymmetric subsonic flow through the nozzle
which satisfies $\eqref{dd2}$ and the properties in {\rm (ii)}.

\item[(iv)]
There exists a critical mass flux $\hat{m}$ such that,
for any $m\in(\delta^{{1}/{4}}, \hat{m})$,
there exists a unique axisymmetric subsonic flow through the nozzle
with the mass flux condition $(\ref{dd-9})$ and the asymptotic
behavior $(\ref{dd4})$. Moreover, $\hat{m}$ is the upper critical mass
flux for the existence of subsonic flow in the following sense:
Either $\sup\limits_{\overline{\Omega}}(U^2+V^2-c^2)\rightarrow 0$ as $m\rightarrow \hat{m}$,
or there is no $\sigma>0$ such that, for all $m\in(\hat{m} ,\hat{m} +\sigma)$,
there is an Euler flow with the mass flux $m$ through the nozzle which
satisfies the upstream asymptotic condition $(\ref{dd-13})$--$(\ref{dd-13s})$,
the downstream asymptotic behavior $(\ref{dd4})$,
and
$
\sup\limits_{m\in(\hat{m},\hat{m}+\sigma)}\sup\limits_{\overline{\Omega}}
(c^2-(U^2+V^2))>0.
$
\end{enumerate}
\end{theorem}

As above, we have the subsonic-sonic limit theorem for this case.

\begin{theorem}[Subsonic-sonic limit of three-dimensional Euler flows through an axisymmetric nozzle]\label{thm4.3}
Let $\delta^{{1}/{4}}<m^{\varepsilon}<\hat{m}$ be a sequence of mass fluxes, and let
$\rho^\varepsilon$, $u^\varepsilon=(u^\varepsilon_1, u^\varepsilon_2, u^\varepsilon_3)$,
and $p^\varepsilon$ be the corresponding solutions to $\mathbf{Problem~2}~(m^\varepsilon)$.
Then, as $m^{\varepsilon}\rightarrow \hat{m}$, the solution
sequence  possesses a subsequence (still denoted by)
$(\rho^\varepsilon, u^{\varepsilon}, p^\varepsilon)$
that converges strongly {\it a.e.} in $\Omega$ to
a vector function
$(\rho, u, p)$ with $u=(u_1, u_2, u_3)$ which is a weak solution of
$\mathbf{Problem~2} (\hat{m})$.
Furthermore, the limit solution $(\rho, u, p)$ also satisfies $(\ref{1.5})$
in the distributional sense and the boundary conditions $(\ref{cdx-sc})$
as the normal trace of the divergence-measure field
$(\rho u_1,\rho u_2, \rho u_3)$ on the boundary in the sense of Chen-Frid {\rm \cite{Chen7}}.
\end{theorem}

\noindent{\bf Proof}.
First, we need to show that
$(\rho^\varepsilon, u^\varepsilon, p^\varepsilon)$ satisfy condition $(A.1)$--$(A.3)$
in $\Omega$.

For the approximate solutions, $B^\varepsilon$ and $S^\varepsilon$
satisfy
\begin{eqnarray}
&& \partial_{x_1}(r\rho^\varepsilon U^\varepsilon B^\varepsilon)+\partial_r (r\rho^\varepsilon V^\varepsilon B^\varepsilon) = 0,\\[2mm]
&& \partial_{x_1}(r\rho^\varepsilon U^\varepsilon S^\varepsilon)+\partial_r (r\rho^\varepsilon V^\varepsilon S^\varepsilon) = 0.
\end{eqnarray}
From $\partial_{x_1}(r\rho^\varepsilon U^\varepsilon)+\partial_r (r\rho^\varepsilon V^\varepsilon)=0$, we
introduce $\psi^\varepsilon$ as
\begin{eqnarray}
\begin{cases}
\partial_{x_1}\psi^\varepsilon=-r \rho^\varepsilon V^\varepsilon,\\[2mm]
\partial_{r}\psi^\varepsilon=r \rho^\varepsilon U^\varepsilon.
\end{cases}
\end{eqnarray}
From the far-field behavior of the Euler flows,
we define $\psi^\varepsilon_-(r):=\lim\limits_{x_1\rightarrow-\infty}\psi^\varepsilon(x_1, r)$.
Similar to the argument in Theorem \ref{thm5.2}, $(\psi^\varepsilon_{-})^{-1}(\psi^\varepsilon)$ are
nondegenerate maps.
A direct calculation yields
\begin{eqnarray*}
&& B^\varepsilon(x_1, x_2, x_3)=B_{-}((\psi^\varepsilon_{-})^{-1}(\psi^\varepsilon(x_1, \sqrt{x_2^2+x_3^2}))),\\[2mm]
&& S^\varepsilon(x_1, x_2, x_3)=S_{-}((\psi^\varepsilon_{-})^{-1}(\psi^\varepsilon(x_1, \sqrt{x_2^2+x_3^2}))).
\end{eqnarray*}
Therefore, we have
\begin{eqnarray}
\begin{cases}
\partial_{x_1}B^\varepsilon=-r\rho^\varepsilon V^\varepsilon\frac{B'_{-}}{(\psi^\varepsilon_-)'},\\[2mm]
\partial_{x_2}B^\varepsilon= x_2\rho^\varepsilon U^\varepsilon\frac{B'_{-}}{(\psi^\varepsilon_-)'},\\[2mm]
\partial_{x_3}B^\varepsilon=x_3\rho^\varepsilon U^\varepsilon\frac{B'_{-}}{(\psi^\varepsilon_-)'}.
\end{cases}
\end{eqnarray}
Notice that
\begin{equation}
\frac{B'_-}{(\psi^\varepsilon_-)'}((\psi^\varepsilon_-)^{-1}(\psi^\varepsilon))
=\frac{B'_-((\psi^\varepsilon_-)^{-1}(\psi^\varepsilon))}
{(\psi^\varepsilon_-)^{-1}(\psi^\varepsilon)\rho^\varepsilon_-U^\varepsilon_-((\psi^\varepsilon_-)^{-1}(\psi^\varepsilon))}.
\end{equation}
Since $B'_-(0)=0$ and $B_-\in C^{1,1}$, we conclude that $\frac{B'_-(s)}{s}$ is bounded.
Then the sequence $B^\varepsilon$ is uniformly bounded in $BV$, which implies its strong convergence.
The similar argument can lead to the strong convergence of $S^\varepsilon$.

\smallskip
On the other hand, the vorticity $\omega^\varepsilon$ has the following expression:
\begin{eqnarray}
\begin{cases}
\omega_{1,2}^\varepsilon=\partial_{x_1}u_2^\varepsilon-\partial_{x_2}u_1^\varepsilon
   =\frac{x_2}{r}(\partial_{x_1}V^\varepsilon-\partial_{r}U^\varepsilon),\\[2mm]
\omega_{2,3}^\varepsilon=\partial_{x_2}u_3^\varepsilon-\partial_{x_3}u_2^\varepsilon=0,\\[2mm]
\omega_{3,1}^\varepsilon=\partial_{x_3}u_1^\varepsilon-\partial_{x_1}u_3^\varepsilon
 =-\frac{x_3}{r}(\partial_{x_1}V^\varepsilon-\partial_{r}U^\varepsilon).
\end{cases}
\end{eqnarray}
A direct computation gives
\begin{equation}
\partial_{x_1}V^\varepsilon-\partial_{r}U^\varepsilon
=\frac{r}{(\psi^\varepsilon_-)'}\big(\rho^\varepsilon B_-'-\frac{(\rho^\varepsilon)^\gamma S_-'}{\gamma}\big),
\end{equation}
which implies that $\omega^\varepsilon$ is uniformly bounded in the bounded measure space.

Since
$(\ref{1.5})$ holds for the sequence of subsonic solutions
$(\rho^\varepsilon, u^{\varepsilon}, p^\varepsilon)(x)$,
it is straightforward to see from Theorem $\ref{thm3.1}$ that there exists
a subsequence (still denoted by) $(\r^\v,u^\v,p^\v)$ which converges to a vector function $(\r,u,p)$ {\it a.e.} in $\Omega$ satisfying $(\ref{1.5})$ in the distributional sense.

\smallskip
The boundary condition is satisfied
in the sense of Chen-Frid \cite{Chen7}, which implies
\begin{equation}
\int_{\partial \Omega}\phi(w) (\rho u)(w)\cdot \nu (w)\, d\mathcal{H}^{1}(w)
= \int_{\Omega} (\rho u)(x) \cdot \nabla \phi(x)\, dx + \langle \mbox{div}(\rho u)|_{\Omega}, \phi\rangle
\end{equation}
for $\psi\in\mathbf{C}^1_0$. From above, we can see
$\langle \mbox{div}(\rho u)|_{\Omega}, \phi\rangle=0$. Furthermore, we have
\begin{equation}
\int_{\Omega}(\rho u)(x) \cdot \nabla \phi(x)\, dx
=\lim\limits_{\varepsilon \rightarrow 0}
\int_{\Omega}(\rho^\varepsilon u^\varepsilon)(x)\cdot \nabla \phi(x)\, dx=0,
\end{equation}
which yields
\begin{equation}
\int_{\partial \Omega}\phi(w)(\rho u)(w)\cdot\nu(w) \, d\mathcal{H}^{1}(w)=0,
\end{equation}
that is, $(\rho u)\cdot \nu = 0$ on $\partial\Omega$ in $\mathcal{D}'$.

This completes the proof.

\begin{remark}
In the homentropic case, the subsonic results of \cite{CX,DD,Xin4} can be also extended to the subsonic-sonic limit by
using Theorem $\ref{thm2.1}$.
\end{remark}

\medskip
\noindent {\bf Acknowledgments:}
The research of
Gui-Qiang Chen was supported in part by
the UK EPSRC Science and Innovation
Award to the Oxford Centre for Nonlinear PDE (EP/E035027/1),
the NSFC under a joint project Grant 10728101, and
the Royal Society--Wolfson Research Merit Award (UK).
The research
of Fei-Min Huang was supported in part by NSFC Grant No. 10825102 for
distinguished youth scholars, and the National Basic Research Program of
China (973 Program) under Grant No. 2011CB808002.
The research of Tian-Yi Wang was supported in part
by the China Scholarship Council  No. 201204910256
as an exchange graduate student at the University of Oxford, the UK EPSRC Science and Innovation
Award to the Oxford Centre for Nonlinear PDE (EP/E035027/1), and the NSF of China under Grant 11371064.

\end{document}